\documentclass{amsart}

\usepackage{amssymb}
\usepackage{amsthm}
\usepackage{amsmath}
\usepackage[all]{xy} \SelectTips{eu}{}
\usepackage[hidelinks]{hyperref}

\newcommand{\numberseries}{\bfseries} 

\newlength{\thmtopspace} 
\newlength{\thmbotspace} 
\newlength{\thmheadspace} 
\newlength{\thmindent} 

\newtheoremstyle{bfupright head,upright body}
{\thmtopspace}{\thmbotspace}
{\upshape}{\thmindent}{\bfseries}{.}{\thmheadspace} {{\numberseries
    \thmnumber{#2\;}}\thmnote{#3}}

\newtheoremstyle{fixed bf head,slanted body}
{\thmtopspace}{\thmbotspace}{\slshape}
{\thmindent}{\bfseries}{.}{\thmheadspace} {{\numberseries
    \thmnumber{#2\;}}\thmname{#1}\thmnote{ (#3)}}

\newtheoremstyle{fixed bf head,upright body}
{\thmtopspace}{\thmbotspace}{\upshape}
{\thmindent}{\bfseries}{.}{\thmheadspace} {{\numberseries
    \thmnumber{#2\;}}\thmname{#1}\thmnote{ (#3)}}

\theoremstyle{bfupright head,upright body} \newtheorem{res}{}[section]

\theoremstyle{fixed bf head,slanted body}
\newtheorem{thm}[res]{Theorem} \newtheorem*{thm*}{Theorem}
\newtheorem{prp}[res]{Proposition} \newtheorem*{prp*}{Proposition}
\newtheorem{cor}[res]{Corollary} \newtheorem*{cor*}{Corollary}
\newtheorem{lem}[res]{Lemma} \newtheorem*{lem*}{Lemma}

\theoremstyle{fixed bf head,upright body}
\newtheorem{dfn}[res]{Definition} \newtheorem*{dfn*}{Definition}
\newtheorem{rmk}[res]{Remark} \newtheorem*{rmk*}{Remark}

\newlength{\thmlistleft} 
\newlength{\thmlistright} 
\newlength{\thmlistpartopsep} 
\newlength{\thmlisttopsep} 
\newlength{\thmlistparsep} 
\newlength{\thmlistitemsep} 

\newcounter{eqc} \newenvironment{eqc}{\begin{list}{\upshape
      (\textit{\roman{eqc}})}%
    {\usecounter{eqc}%
      \setlength{\leftmargin}{\thmlistleft}%
      \setlength{\labelwidth}{\thmlistleft}%
      \setlength{\rightmargin}{\thmlistright}%
      \setlength{\partopsep}{\thmlistpartopsep}%
      \setlength{\topsep}{\thmlisttopsep}%
      \setlength{\parsep}{\thmlistparsep}%
      \setlength{\itemsep}{\thmlistitemsep}}}%
  {\end{list}}%

\newcommand{\eqclbl}[1]{{\upshape(\textit{#1})}}

\newenvironment{prf*}[1][Proof]{%
  \begin{proof}[\bf #1]
    \setcounter{equation}{0}
    } {\end{proof} }

\newcommand{\proofofimp}[3][:]{\mbox{\eqclbl{#2}$\!\implies\!$\eqclbl{#3}#1}}

\newcommand{\thmref}[2][Theorem~]{#1\ref{thm:#2}}
\newcommand{\corref}[2][Corollary~]{#1\ref{cor:#2}}

\newcommand{\lemref}[2][Lemma~]{#1\ref{lem:#2}}

\newcommand{\rmkref}[2][Remark~]{#1\ref{rmk:#2}}
\newcommand{\secref}[2][Section~]{#1\ref{sec:#2}}

\newcommand{\thmcite}[2][?]{\cite[Thm.~#1]{#2}}
\newcommand{\corcite}[2][?]{\cite[Cor.~#1]{#2}}
\newcommand{\prpcite}[2][?]{\cite[Prop.~#1]{#2}}
\newcommand{\lemcite}[2][?]{\cite[Lem.~#1]{#2}}

\newcommand{\setof}[3][\mspace{1mu}]{\{#1#2 \mid #3#1\}}
\newcommand{\ZZ}{\mathbb{Z}}
\newcommand{\deq}{\:=\:}
\newcommand{\dle}{\:\le\:}
\newcommand{\dis}{\:\is\:}
\renewcommand{\d}{v}            
\newcommand{\is}{\cong}
\newcommand{\qis}{\simeq}
\renewcommand{\le}{\leqslant}
\renewcommand{\ge}{\geqslant}
\newcommand{\lra}{\longrightarrow}
\newcommand{\xra}[2][]{\xrightarrow[#1]{\;#2\;}}
\newcommand{\qra}{\xra{\qis}}
\newcommand{\Aop}{A^\circ}
\newcommand{\dif}[2][]{{\partial}^{#2}_{#1}}
\newcommand{\Cy}[2][]{\operatorname{Z}_{#1}(#2)}
\newcommand{\Co}[2][]{\operatorname{C}_{#1}(#2)}
\renewcommand{\H}[2][]{\operatorname{H}_{#1}(#2)}
\newcommand{\Tha}[2]{#2_{{\scriptscriptstyle\le}#1}}
\newcommand{\Thb}[2]{#2_{{\scriptscriptstyle\ge}#1}}
\newcommand{\fd}[2][A]{\operatorname{fd}_{#1}#2}
\newcommand{\id}[2][A]{\operatorname{id}_{#1}#2}
\newcommand{\pd}[2][A]{\operatorname{pd}_{#1}#2}
\newcommand{\Gfd}[2][A]{\operatorname{Gfd}_{#1}#2}
\newcommand{\Gpd}[2][A]{\operatorname{Gpd}_{#1}#2}
\newcommand{\Gfcd}[2][A]{\operatorname{Gfcd}_{#1}#2}
\newcommand{\Ggldim}[1][A]{\operatorname{Ggldim}(#1)}
\newcommand{\Gwgldim}[1][A]{\operatorname{Gwgldim}(#1)}
\newcommand{\FPD}[1][A]{\operatorname{FPD}(#1)}
\newcommand{\silp}[1][A]{\operatorname{silp}(#1)}
\newcommand{\silf}[1][A]{\operatorname{silf}(#1)}
\newcommand{\spli}[1][A]{\operatorname{spli}(#1)}
\newcommand{\splf}[1][A]{\operatorname{splf}(#1)}
\newcommand{\sfli}[1][A]{\operatorname{sfli}(#1)}
\newcommand{\Hom}[3][A]{\operatorname{Hom}_{#1}(#2,#3)}
\newcommand{\Ext}[4][A]{\operatorname{Ext}_{#1}^{#2}(#3,#4)}
\def\urltilda{\kern -.15em\lower .7ex\hbox{\~{}}\kern .04em}

\title{Gorenstein weak global dimension is symmetric}

\author[L.W.\ Christensen]{Lars Winther Christensen} \address{L.W.C. \
  Texas Tech University, Lubbock, TX 79409, U.S.A.}
\email{lars.w.christensen@ttu.edu}
\urladdr{http://www.math.ttu.edu/\urltilda lchriste}

\author[S.\ Estrada]{Sergio Estrada} \address{S.E. \ Universidad de
  Murcia, Murcia 30100, Spain} \email{sestrada@um.es}
\urladdr{https://webs.um.es/sestrada/}

\author[P.\ Thompson]{Peder Thompson} \address{P.T. \ Norwegian
  University of Science and Technology, 7491 Trondheim, Norway}
\email{peder.thompson@ntnu.no} \urladdr{https://folk.ntnu.no/pedertho}

\thanks{L.W.C.\ was partly supported by Simons Foundation
  collaboration grant 428308. S.E.\ was partly supported by grant
  MTM2016-77445-P and FEDER funds and by the grant 19880/GERM/15 from
  the Fundaci\'on S\'eneca-Agencia de Ciencia y Tecnolog\'{\i}a de la
  Regi\'on de Murcia.}

\date{3 June 2021}



\begin{document}

\begin{abstract}
  We study the Gorenstein weak global dimension of associative rings
  and its relation to the Gorenstein global dimension. In particular,
  we prove the conjecture that the Gorenstein weak global dimension is
  a left--right symmetric invariant---just like the (absolute) weak
  global dimension.
\end{abstract}

\maketitle

\section*{Introduction}

\noindent

\noindent
A guiding principle in Gorenstein homological algebra is to seek
analogues of results about absolute homological dimensions. For
example, the Gorenstein global dimension of a ring can equally well be
computed in terms of the Gorenstein projective or Gorenstein injective
dimensions of its modules; this was proved by Enochs and Jenda
\cite{rha} in the noetherian case and by Bennis and Mahdou
\cite{DBnNMh10} in general.

The notion of a weak global dimension has also been considered in
Gorenstein homological algebra, for example by Emmanouil
\cite{IEm12}. The weak global dimension of a ring is a left--right
symmetric invariant; that is, a ring has finite weak global dimension
on the left if and only if it enjoys the same property on the
right. In \secref{1} we prove the corresponding statement in
Gorenstein homological algebra, thus confirming a widely held
conjecture that was formally stated by Bennis \cite{DBn10}.

In \secref{2} we use this symmetry to investigate the relations
between the Gorenstein global and Gorenstein weak global
dimensions. The main result of this section, \thmref{main}, shows that
finite Gorenstein weak global dimension together with finite
projective dimension of flat modules implies finite Gorenstein global
dimension; under extra assumptions on the ring this was proved by
Bennis and Mahdou~\cite{DBnNMh09}. The same theorem relates the
Gorenstein global dimension to a new invariant: the Gorenstein
flat-cotorsion dimension, which was introduced in \cite{CELTWY-21}. In
fact, this new invariant plays a key role already in the proof of
\thmref{sym}, the main result of the first section. The invariant is
built on the theory of Gorenstein flat-cotorsion modules developed in
\cite{CET-20} as well as recent work of \v{S}aroch and
\v{S}tov{\'{\i}}{\v{c}}ek \cite{JSrJSt20}.
\begin{equation*}
  \ast \ \ \ast \ \ \ast
\end{equation*}
Throughout the paper, $A$ denotes an associative ring. By an
$A$-module we mean a left $A$-module, and we treat right $A$-modules
as modules over the opposite ring~$\Aop$. By an \emph{$A$-complex} we
mean a complex of $A$-modules. For such a complex $M$ and an integer
$n$, the hard truncation of $M$ above at $n$ is denoted $\Tha{n}{M}$
while $\Thb{n}{M}$ denotes the hard truncation of $M$ below at $n$.
For $v \in \ZZ$ the cycle module in degree $v$, i.e.\ the kernel of
$\dif[v]{M}$ is denoted $\Cy[v]{M}$, while $\Co[v]{M}$ denotes the
cokernel module in degree $v$, i.e.\ the cokernel of
$\dif[v+1]{M}$. We say that $M$ has \emph{bounded} homology if
$\H[v]{M}=0$ holds for $|v| \gg 0$.

The notation and terminology above is all standard; the only
non-standard terminology applied in this paper comes from
\cite{CELTWY-21,CET-20}: An acyclic complex $T$ of flat-cotorsion
$A$-modules is called \emph{totally acyclic} if $\Hom{T}{F}$ is
acyclic for every flat-cotorsion $A$-module $F$. A cycle in such a
complex is called a \emph{Gorenstein flat-cotorsion} module.  A
semi-flat complex of flat-cotorsion modules is called
\emph{semi-flat-cotorsion.} In the derived category of $A$, every
complex $M$ is isomorphic to a semi-flat-cotorsion complex---see for
example \cite[Cnstr.~2.4]{CELTWY-21}; such a complex is called a
semi-flat-cotorsion replacement of $M$. The Gorenstein flat-cotorsion
dimension of an $A$-complex $M$ is denoted $\Gfcd{M}$; it is the least
$n \ge \sup\setof{v\in\ZZ}{\H[v]{M} \ne 0}$ such that the
$n^{\mathrm{th}}$ cokernel in a semi-flat-cotorsion replacement of $M$
is Gorenstein flat-cotorsion.

For the absolute homological dimensions we use two-letter
abbreviations---pd, id, and fd---and we write Gpd, Gid, and Gfd for
the corresponding Gorenstein dimensions. The notation for the
invariant
\begin{equation*}
  \splf \deq \sup\setof{\pd{F}}{F \text{ is a flat $A$-module}}
\end{equation*}
is another acronym, ``splf'' stands for ``supremum of projective
lengths of flat modules.'' The invariants sfli, spli, silp, and silf
are defined similarly; see \cite[\S 1.2]{IEm12}.

\section{Symmetry of Gorenstein weak global dimension}
\label{sec:1}

\noindent
Holm proves in \thmcite[2.6]{HHl04c} that if $A$ is coherent and
$\splf[\Aop]$ is finite, then the equality $\Gfd{M} = \fd{M}$ holds
for $A$-modules of finite injective dimension. The key to our proof of
the main result in this section is to show that this equality holds
without the assumptions on $A$. By the work done in \cite{CELTWY-21}
it suffices to prove the analogous equality for the Gorenstein
flat-cotorsion dimension, and since this is a result of independent
interest, we prove it for complexes.

\begin{thm}
  \label{thm:Gfd-fd-for-id-finite}
  Let $M$ be an $A$-complex with bounded homology. If $M$ has finite
  injective dimension, then the equality $\Gfcd{M}=\fd{M}$ holds.
\end{thm}

\begin{prf*}
  The equality $\Gfcd{M} = \fd{M}$ holds trivially if $M$ is acyclic,
  so assume that $M$ is not acyclic and assume further, without loss
  of generality, that $\id{M}=0$ holds.  Set $w = \sup\setof{v\in\ZZ}{\H[v]{M} \ne 0}$ and let
  $M \qra I$ be a semi-injective resolution with $I_\d = 0$ for
  $\d > w$ and $\d < 0$.  There is an exact sequence of complexes
  $0\to C' \to F \to I \to 0$ with $F$ semi-flat-cotorsion and $C'$ an
  acyclic complex of cotorsion modules; this follows from work of
  Gillespie \cite{JGl04}, see also \cite[Fact 2.2]{CELTWY-21}.
  Acyclicity of $C'$ yields an exact sequence
  $0\to \Cy[0]{C'}\to \Cy[0]{F}\to I_0\to 0$. Since both $\Cy[0]{C'}$
  and $I_0$ are cotorsion---for the former see Bazzoni,
  Cort\'es-Izurdiaga, and Estrada \thmcite[1.3]{BCE-20}---so is
  $\Cy[0]{F}$.  There is thus a semi-flat-cotorsion resolution
  $F'\to \Cy[0]{F}$ concentrated in non-negative degrees, constructed
  by taking successive flat covers. This complex glued together with
  $\Tha{0}{F}$ is acyclic and semi-flat, see Christensen and Holm \cite[6.1]{LWCHHl15}, so
  per \thmcite[7.3]{LWCHHl15} the module $\Cy[0]{F}$ is
  flat-cotorsion, and we may assume that $F_\d=0$ holds for $\d<0$.

  Fix $n\ge \Gfcd{M}$; the module $\Co[n]{F}$ is Gorenstein
  flat-cotorsion by \lemcite[4.3]{CELTWY-21}. We argue next that
  $\Ext{1}{G}{\Co[n]{F}} = 0$ holds for every Gorenstein
  flat-cotorsion module $G$. Fix such a module $G$. By definition,
  there is an exact sequence of $A$-modules,
  $0\to G \to T_0 \to \cdots \to T_{-n+1} \to G' \to 0$, with each
  $T_\d$ flat-cotorsion and $G'$ Gorenstein flat-cotorsion. As
  $\Co[n]{F}$ is cotorsion, dimension shifting yields:
  \begin{equation*}
    \Ext{1}{G}{\Co[n]{F}} \dis \Ext{n+1}{G’}{\Co[n]{F}} \:.
  \end{equation*}
  Let $C$ be the mapping cone of the quasi-isomorphism $F\to I$; it is
  concentrated in non-negative degrees and consists of sums of
  modules that are flat-cotorsion or injective. Moreover, one has
  $C_0 = I_0$ as $F_v = 0$ for $v<0$, and because $I_v=0$ for
  $v > w$ and $n \ge w$ holds, there is an isomorphism
  $\Co[n+1]{C} \is \Co[n]{F}$.  Thus dimension shifting along
  \begin{equation*}
    0 \lra \Co[n+1]{C} \lra C_{n} \lra \cdots \lra C_{1} \lra I_0 \lra 0
  \end{equation*}
  yields
  \begin{equation*}
    \Ext{n+1}{G’}{\Co[n]{F}} \dis \Ext{1}{G’}{I_0} \deq 0 \:.
  \end{equation*}
  Combining the displayed isomorphisms one gets
  $\Ext{1}{G}{\Co[n]{F}}=0$. Now, with $g = \Gfcd{M}$ and $n=g+1$ and
  $G = \Co[g]{F}$ one has $\Ext{1}{\Co[g]{F}}{\Co[g+1]{F}}=0$. This
  means that the exact sequence
  $0 \to \Co[g+1]{F} \to F_g \to \Co[g]{F} \to 0$ splits, whence
  $\Co[g]{F}$ is flat-cotorsion. Thus one has $\fd{M} \le g$, and the
  opposite inequality holds by \lemcite[5.11]{CELTWY-21}.
\end{prf*}

In particular we now have the desired strengthening of
\thmcite[2.6]{HHl04c}.

\begin{cor}
  \label{cor:Gfd-fd-for-id-finite}
  Let $M$ be an $A$-complex with bounded homology. If $M$ has finite
  injective dimension, then the equality $\Gfd{M}=\fd{M}$ holds.
\end{cor}

\begin{prf*}
  The Gorenstein flat dimension is a refinement of the flat dimension,
  so if $\Gfd{M}=\infty$ holds, then the equality is trivial. If
  $\Gfd{M}<\infty$, then \thmcite[5.7]{CELTWY-21} yields
  $\Gfcd{M}=\Gfd{M}$ and the asserted equality follows from
  \thmref{Gfd-fd-for-id-finite}.
\end{prf*}

The Gorenstein global dimension of $A$, denoted $\Ggldim$, is the
supremum of the Gorenstein projective dimensions (equivalently, see
\thmcite[1.1]{DBnNMh10}, the Gorenstein injective dimensions) of all
$A$-modules.

\begin{dfn}
The
\emph{Gorenstein weak global dimension}~of~$A$ is
\begin{equation*}
  \Gwgldim \deq \sup\setof{\Gfd{M}}{M \text{ is an $A$-module}} \:.
\end{equation*}
This is the invariant that Bennis and Mahdou denote
$l.\operatorname{wGgldim}(A)$ in \cite{DBn10,DBnNMh10} and
$\operatorname{G-wdim}(A)$ in \cite{DBnNMh09}. When $\Gwgldim$ and
$\Gwgldim[\Aop]$ are finite, and only then, Emmanouil \cite{IEm12}
uses the symbol $\operatorname{Gw.dim}A$ for their common value.
\end{dfn}

If $\Gwgldim$ is finite then so is $\sfli[\Aop]$; this is elementary,
see \lemcite[5.1]{IEm12}. On the other hand, if both $\sfli$ and
$\sfli[\Aop]$ are finite, then per \thmcite[5.3]{IEm12} both
$\Gwgldim$ and $\Gwgldim[\Aop]$ are finite. Thus, the key to prove
symmetry of the Gorenstein weak global dimension is to see that
$\Gwgldim < \infty$ implies $\sfli < \infty$. In our proof of
\thmref{sym} this follows from \corref{Gfd-fd-for-id-finite}, which
through \thmcite[5.7]{CELTWY-21} relies crucially on the work of
\v{S}aroch and \v{S}tov{\'{\i}}{\v{c}}ek \cite{JSrJSt20}. In
\rmkref{SS} we sketch how to obtain
symmetry directly from \cite{JSrJSt20}. However, there is more to
\thmref{sym}: In the next section it facilitates the comparison of
$\Gwgldim$ to $\Ggldim$, see for example \corref{comm}.

\begin{thm}
  \label{thm:sym}
  The following conditions are equivalent.
  \begin{eqc}
  \item $\Gwgldim < \infty$.
  \item All $A$-modules have finite Gorenstein flat dimension.
  \item $\sfli$ and $\sfli[\Aop]$ are finite.    
  \item All $A$- and $\Aop$-modules have finite Gorenstein flat
    dimension.
  \item All $A$- and $\Aop$-modules have finite Gorenstein
    flat-cotorsion dimension.
  \end{eqc}
\end{thm}

\begin{prf*}
  Conditions $(i)$ and $(ii)$ are equivalent as the class of Gorenstein
  flat modules is closed under coproducts; see
  for example Holm \prpcite[3.2]{HHl04a}. Conditions $(iii)$ and $(iv)$ are equivalent
  by \thmcite[5.3]{IEm12}, and $(iv)$ evidently implies $(ii)$.

  \proofofimp{ii}{iii} It follows from
  \corref{Gfd-fd-for-id-finite} that $\sfli$ is finite and from
  \lemcite[5.1]{IEm12} that $\sfli[\Aop]$ is finite.

  \proofofimp{iv}{v} The Gorenstein flat-cotorsion dimension of a
  module is bounded above by its Gorenstein flat dimension, see
  \thmcite[5.7]{CELTWY-21}.

  \proofofimp{v}{iii} It follows from
  \thmref{Gfd-fd-for-id-finite} that $\sfli$ and $\sfli[\Aop]$ are
  finite.
\end{prf*}

The next equality was formally conjectured by Bennis
\cite[Conj.~1.1]{DBn10}.  For emphasis, we point out that it shows
that Emmanouil's \cite{IEm12} definition can be relaxed: It suffices
to consider finiteness of Gorenstein flat dimensions on one side of
the ring.

\begin{cor}
  \label{cor:sym}
  One has $\Gwgldim = \Gwgldim[\Aop]$.
\end{cor}

\begin{prf*}
  The invariants $\Gwgldim$ and $\Gwgldim[\Aop]$ are simultaneously
  finite by \thmref{sym}, and when finite they are equal by
  \thmcite[5.3]{IEm12}.
\end{prf*}

\begin{rmk}
  \label{rmk:SS}
  That $\Gwgldim < \infty$ implies $\sfli < \infty$ can be deduced
  directly from \cite{JSrJSt20}: In the notation of that paper, given
  a module $M\in \mathcal{PGF}^\perp$, there exists by
  \thmcite[4.9]{JSrJSt20} an exact sequence,
  $0\to H\to T_{n-1}\to \cdots \to T_0\to M\to 0$, with each $T_v$ a
  projective $A$-module and $H$ in $\mathcal{PGF}^\perp$. If
  $\Gfd{M}\le n$, then $H$ is also Gorenstein flat, hence flat per
  \thmcite[4.11]{JSrJSt20}.
\end{rmk}

\section{Comparing Gorenstein global dimensions}
\label{sec:2}

\noindent
In this section, we consider relations between finiteness of the
Gorenstein global dimensions. We begin with a key lemma that compares
the relevant invariants at the level of (complexes of) modules.

\begin{lem}
  \label{lem:Gfd-Gpd}
  For every $A$-complex $M$ with $\H{M} \ne 0$ one has
  \begin{equation*}
    \Gpd{M} \dle \Gfd{M} + \splf \:.
  \end{equation*}

\end{lem}

\begin{prf*}
  Set $n = \splf$ and assume that it is finite. Let $M$ be an
  $A$-complex with $\Gfd{M} = d$ for some integer $d$. Let $P \to M$
  be a semi-projective resolution; the module $C = \Co[d]{P}$ is
  Gorenstein flat---see Christensen, K\"oksal, and Liang
  \prpcite[5.12]{CKL-17}\footnote{Every ring is GF-closed by
    \corcite[4.12]{JSrJSt20}.}---and it suffices to show that
  $\Gpd{C} \le n$ holds, as this implies that $\Co[d+n]{P}$ is
  Gorenstein projective. By assumption there is an acyclic complex,
  $0\to C \to F_0\to F_{-1} \to \cdots$, with each module $F_\d$ flat
  and each cokernel Gorenstein flat. As in Cartan and Eilenberg's
  \cite[Chapter XVII, \S1]{CarEil}, or the proof of
  \lemcite[5.2]{IEm12}, construct a projective resolution of this
  complex in the category of $A$-complexes:
  \begin{equation*}
    \xymatrix@=1.5pc{
      & \vdots\ar[d] & \vdots\ar[d] & \vdots\ar[d] &  \\
      0 \ar[r] & Q_1 \ar[r]\ar[d] & Q_1^{(0)} \ar[r] \ar[d]& Q_1^{(-1)} \ar[r]\ar[d] & \cdots\\
      0 \ar[r] & Q_0 \ar[r]\ar[d] & Q_0^{(0)} \ar[r]\ar[d] & Q_0^{(-1)} \ar[r]\ar[d] & \cdots\\
      0 \ar[r] &C \ar[r] \ar[d] & F_0 \ar[r] \ar[d] & F_{-1} \ar[r] \ar[d] & \cdots\\
      & 0 & 0 & 0 & 
    }
  \end{equation*}
  This induces an exact sequence
  \begin{equation*}
    0 \lra \Co[n]{Q} \lra \Co[n]{Q^{(0)}} \lra \Co[n]{Q^{(-1)}} \lra \cdots \:.
  \end{equation*}
  The class of Gorenstein flat $A$-modules is resolving by
  \corcite[4.12]{JSrJSt20}, so the module $\Co[n]{Q}$ is Gorenstein
  flat. By assumption, $\Co[n]{Q^{(i)}}$ is projective for $i\le 0$,
  and by construction the cokernels of the exact sequence are
  Gorenstein flat, so $\Co[n]{Q}$ is Gorenstein projective by
  \thmcite[4.4]{JSrJSt20}. Thus $\Gpd{C} \le n$ holds as desired.
\end{prf*}

Jiangsheng Hu pointed us to the following easy consequence of
\lemref{Gfd-Gpd}.

\begin{prp}
  Let $n$ be an integer. The following conditions are equivalent.
  \begin{eqc}
  \item $\splf \le n$.
  \item Every Gorenstein flat $A$-module has Gorenstein projective
    dimension at most $n$ and every flat Gorenstein projective
    $A$-module is projective.
  \item Every flat $A$-module has Gorenstein projective dimension at
    most $n$ and every flat Gorenstein projective $A$-module is
    projective.
  \end{eqc}
\end{prp}

\begin{prf*}
  Evidently, $(ii)$ implies $(iii)$.

  \proofofimp{i}{ii} For a Gorenstein flat $A$-module $M$,
  \lemref{Gfd-Gpd} yields $\Gpd{M} \le n$. A flat Gorenstein
  projective $A$-module $G$ is projective as $\Gpd{G} = \pd{G}$ holds
  because the Gorenstein projective dimension refines the projective
  dimension.

  \proofofimp{iii}{i} The $n^\mathrm{th}$ syzygy of a flat $A$-module
  is Gorenstein projective and flat, hence projective.
\end{prf*}

This brings us to the main result of this section; it compares to
\thmcite[2.1]{DBnNMh09} as does \corref{coh_cor}.

\begin{thm}
  \label{thm:main}
  There are inequalities
  \begin{equation*}
    \sup\setof{\Gfcd{M}}{M \text{ is an $A$-module}}
    \dle \Ggldim \dle \Gwgldim + \splf \:.
  \end{equation*}
\end{thm}

\begin{prf*}
  The second inequality follows immediately from \lemref{Gfd-Gpd}. To
  prove the first inequality, set $n=\Ggldim$; we may assume that it
  is finite.  By \thmcite[4.1]{IEm12} one has $\spli=n=\silp$, and a
  result of Emmanouil and Talelli~\prpcite[2.1]{IEmOTl11} yields
  $\silf = n$. We first show that every cotorsion $A$-module $C$ has
  $\Gfcd{C}\le n$. To see this, let $C\to I$ be an injective
  resolution. For every $i \le 0$ there is an exact sequence
  $0 \to \Cy[i]{I} \to I_{i} \to \Cy[i-1]{I} \to 0$ of cotorsion
  modules.  Construct flat resolutions $G \to C=\Cy[0]{I}$ and
  $G^{(i)} \to \Cy[i]{I}$ for $i < 0$ by taking successive flat
  covers. By \lemcite[8.2.1]{rha} there are flat resolutions
  $F^{(i)} \to I_i$ which fit into exact sequences
  $0\to G^{(i)} \to F^{(i)} \to G^{(i-1)}\to 0$, such that each module
  $F^{(i)}_j$ is flat-cotorsion, and each syzygy module
  $\Co[j]{F^{(i)}}$ is cotorsion. A standard construction, as in
  \cite[Chap.~ XVII.\S1]{CarEil}, yields a commutative diagram
  \begin{equation*}
    \xymatrix@=1.5pc{
      & \vdots\ar[d] & \vdots\ar[d] & \vdots\ar[d] &  \\
      0 \ar[r] & G_1 \ar[r]\ar[d] & F_1^{(0)} \ar[r] \ar[d]& F_1^{(-1)} \ar[r]\ar[d] & \cdots\\
      0 \ar[r] & G_0 \ar[r]\ar[d] & F_0^{(0)} \ar[r]\ar[d] & F_0^{(-1)} \ar[r]\ar[d] & \cdots\\
      0 \ar[r] & C \ar[r] \ar[d] & I_0 \ar[r] \ar[d] & I_{-1} \ar[r]\ar[d]  & \cdots\\
      & 0 & 0 & 0 &
    }
  \end{equation*}
  with exact rows and columns.  This induces an exact sequence
  \begin{equation*}
    0 \lra \Co[n]{G} \lra \Co[n]{F^{(0)}} \lra \Co[n]{F^{(-1)}} \lra \cdots \:.
  \end{equation*}
  Since $\sfli\le \spli=n$ holds, the modules $\Co[n]{F^{(i)}}$ are
  flat for all $i \le 0$, and by construction they are cotorsion. The
  complex $\Thb{n}{G}$ is a resolution of $\Co[n]{G}$ by
  flat-cotorsion modules, so $\Co[n]{G}$ is a syzygy module in an
  acyclic complex of flat-cotorsion $A$-modules. As $\silf$ is finite,
  this complex is totally acyclic---indeed, for an acyclic complex $X$
  of flat modules and a cotorsion module $Y$ of finite injective
  dimension, a standard dimension shifting argument, as in the proof
  of \thmref{Gfd-fd-for-id-finite}, shows that the complex
  $\Hom{X}{Y}$ is acyclic---so the module $\Co[n]{G}$ is Gorenstein
  flat-cotorsion. Thus $\Gfcd{C}\le n$.
  
  Finally, let $M$ be an $A$-module and $F\to M$ a flat resolution
  built from flat covers. Since $\Co[1]{F}$ is cotorsion, the module
  $\Co[n+1]{F}$ is Gorenstein flat-cotorsion. Work of Gillespie
  \cite{JGl04}, see also \cite[Fact 2.2]{CELTWY-21}, yields an exact
  sequence,
  \begin{equation*}
    0 \lra F \lra C \lra P \lra 0 \:,
  \end{equation*}
  with $C$ degreewise cotorsion and $P$ an acyclic complex of flat
  modules with flat cycle modules. It follows that $C$ is a semi-flat-cotorsion replacement of
  $M$, see \cite[Fact 1.4]{CELTWY-21}.  For $i\ge 1$ the exact
  sequence \mbox{$0 \to F_i\to C_i \to P_i\to 0$} shows that $P_i$ is
  cotorsion. As all the modules $\Co[i]{P}$ are flat, \cite[Lemma
  5.6]{CELTWY-21} applied to $\Thb{1}{P}$ shows that $\Co[n+1]{P}$ is
  flat-cotorsion. It follows that the exact sequence
  \begin{equation*}
    0 \lra \Co[n+1]{F} \lra \Co[n+1]{C} \lra \Co[n+1]{P} \lra 0
  \end{equation*}
  splits, whence $\Co[n+1]{C}$ is Gorenstein flat-cotorsion; in
  particular, $\Gfcd{M}$ is finite.  There is an exact sequence of
  $A$-modules, $0\to M \to C' \to F' \to 0$, with $C'$ cotorsion and
  $F'$ flat. As we have shown above that $\Gfcd{C'}\le n$ holds, it
  now follows from \thmcite[4.5]{CELTWY-21} that also $\Gfcd{M}\le n$.
\end{prf*}

\begin{rmk}
  \label{rmk:FPD}
  Recall that the invariant
  \begin{equation*}
    \FPD \deq \sup\setof{\pd{M}}{M \text{ has finite projective dimension}}
  \end{equation*}
  is known as the \emph{finitistic projective dimension} of $A$. By a
  result of Jensen \prpcite[6]{CUJ70} one has $\splf \le \FPD$.  By
  \thmcite[2.28]{HHl04a} there is an inequality $\FPD \le \Ggldim$,
  and equality holds if $\Ggldim$ is finite.
\end{rmk}

\begin{cor}
  \label{cor:coh_cor}
  Assume that $A$ is right coherent.  There are inequalities
  \begin{equation*}
    \Gwgldim
    \dle \Ggldim \dle \Gwgldim  + \splf \:.
  \end{equation*}
  Moreover, the following conditions are equivalent.
  \begin{eqc}
  \item $\Gwgldim$ and $\splf$ are finite.
  \item $\Ggldim$ is finite.
  \end{eqc}
\end{cor}

\begin{prf*}
  Since $A$ is right coherent, the equality $\Gfcd{M}= \Gfd{M}$ holds
  for every $A$-module $M$ by \corcite[5.8]{CELTWY-21}. Hence, one has
  \begin{equation*}
    \sup\setof{\Gfcd{M}}{M \text{ is an $R$-module}}
    \deq \Gwgldim \;,
  \end{equation*}
  and the asserted inequalities follow from \thmref{main}.

  It is immediate from the second inequality that \eqclbl{i} implies
  \eqclbl{ii}. For the converse, assume that $\Ggldim$ is finite; it
  follows from the first inequality that $\Gwgldim$ is finite, and
  $\splf$ is finite by \rmkref{FPD}.
\end{prf*}

The next corollary applies, in particular, to commutative rings.

\begin{cor}
  \label{cor:comm}
  If $A$ and $\Aop$ are isomorphic, then the next conditions are~\mbox{equivalent.}
  \begin{eqc}
  \item $\Gwgldim$ and $\splf$ are finite.
  \item $\Ggldim$ is finite.
  \end{eqc}
\end{cor}

\begin{prf*}  
  Immediate from \thmref[Theorems~]{sym} and \thmref[]{main}, along
  with \rmkref{FPD}.
\end{prf*}

\begin{cor}
  \label{cor:GGFPD}
  If $\Gwgldim$ is finite, then $\Ggldim = \FPD$ holds, and the
  invariants $\splf$ and $\FPD$ are simultaneously finite.
\end{cor}

\begin{prf*}
  By \rmkref{FPD} the equality holds if $\Ggldim$ is finite. Assume
  that $\Gwgldim$ is finite. If $\FPD$ is finite, then $\splf$ is
  finite by \rmkref{FPD}, so $\Ggldim$ is finite by
  \thmref{main}. This proves the equality, and the last assertion
  follows as finiteness of $\splf$ by \thmref[]{main} implies
  finiteness of $\Ggldim$.
\end{prf*}

Simson \cite{DSm74} shows that a ring $A$ of cardinality $\aleph_n$
has $\splf \le n+1$; for commutative rings this was shown earlier by
Jensen \thmcite[5.8]{CUJ72}. Thus, the next corollary yields, in
particular, that the Gorenstein global dimension is symmetric for
countable coherent rings.

\begin{cor}
  \label{cor:GorSym}
  If $A$ is coherent, then the following conditions are equivalent.
  \begin{eqc}
  \item $\Ggldim$ and $\splf[\Aop]$ are finite.
  \item $\Ggldim[\Aop]$ and $\splf$ are finite.
  \end{eqc}
\end{cor}

\begin{prf*}
  The assertion follows from \corref{coh_cor} combined with
  \corref{sym}.
\end{prf*}

\begin{rmk}
  \label{rmk:1}
  For noetherian rings Beligiannis \corcite[6.11]{ABl00} proved that
  the Gorenstein global dimension is symmetric. \corref{GorSym}
  provides us with non-trivial new examples of rings that exhibit this
  kind of Gorenstein symmetry; by non-trivial we here mean rings of
  infinite global dimension. For instance, let $R$ be a
  non-commutative artinian ring and $S$ a commutative coherent, but
  not noetherian, ring of cardinality $\le \aleph_n$. By \cite{DSm74}
  the direct product ring $R\times S$ satisfies the assumptions in
  \corref{GorSym}. Coherent quotients of the product ring
  $M_n(\mathbb{Q})\times\mathbb{Q}[x_0,x_1,\ldots]$ are simple
  examples of such rings.
\end{rmk}

One can replace the coherent assumption in \corref{GorSym} with
assumptions of Gorenstein flatness of Gorenstein projective modules.

\begin{cor}
  The next conditions are equivalent.
  \begin{eqc}
  \item $\Ggldim$ and $\splf[\Aop]$ are finite and every Gorenstein
    projective $A$-module is Gorenstein flat.
  \item $\Ggldim[\Aop]$ and $\splf$ are finite and every Gorenstein
    projective $\Aop$-module is Gorenstein flat.
  \end{eqc}
\end{cor}

\begin{prf*}
  If $\Ggldim$ is finite and every Gorenstein projective $A$-module is
  Gorenstein flat, then it follows that $\Gwgldim$ is finite. By
  \corref{sym} this implies that $\Gwgldim[\Aop]$ is finite. It
  follows that every cycle in an acyclic complex of flat
  $\Aop$-modules is Gorenstein flat; in particular, every Gorenstein
  projective $\Aop$-module is Gorenstein flat. Finally, \thmref{main}
  yields $\Ggldim[\Aop] < \infty$ and $\splf$ is finite by
  \rmkref{FPD},
\end{prf*}

\begin{rmk}
  The result of Beligiannis, \corcite[6.11]{ABl00}, mentioned in
  \rmkref{1} shows that a noetherian ring is
  Iwanaga--Gorenstein---that is, of finite self-injective dimension on
  both sides---if it has finite Gorenstein global dimension on one
  side. It follows from \thmref{sym} and \corcite[5.10]{CELTWY-21}
  that a noetherian ring is Iwanaga--Gorenstein if it has finite
  Gorenstein weak global dimension on one side. This improves
  \thmcite[12.3.1]{rha}.
\end{rmk}

\section*{Acknowledgment}

\noindent
We thank the referee for pertinent suggestions that improved the
exposition.

\def\soft#1{\leavevmode\setbox0=\hbox{h}\dimen7=\ht0\advance \dimen7
  by-1ex\relax\if t#1\relax\rlap{\raise.6\dimen7
  \hbox{\kern.3ex\char'47}}#1\relax\else\if T#1\relax
  \rlap{\raise.5\dimen7\hbox{\kern1.3ex\char'47}}#1\relax \else\if
  d#1\relax\rlap{\raise.5\dimen7\hbox{\kern.9ex \char'47}}#1\relax\else\if
  D#1\relax\rlap{\raise.5\dimen7 \hbox{\kern1.4ex\char'47}}#1\relax\else\if
  l#1\relax \rlap{\raise.5\dimen7\hbox{\kern.4ex\char'47}}#1\relax \else\if
  L#1\relax\rlap{\raise.5\dimen7\hbox{\kern.7ex
  \char'47}}#1\relax\else\message{accent \string\soft \space #1 not
  defined!}#1\relax\fi\fi\fi\fi\fi\fi}
  \providecommand{\MR}[1]{\mbox{\href{http://www.ams.org/mathscinet-getitem?mr=#1}{#1}}}
  \renewcommand{\MR}[1]{\mbox{\href{http://www.ams.org/mathscinet-getitem?mr=#1}{#1}}}
  \providecommand{\arxiv}[2][AC]{\mbox{\href{http://arxiv.org/abs/#2}{\sf
  arXiv:#2 [math.#1]}}} \def\cprime{$'$}
\providecommand{\bysame}{\leavevmode\hbox to3em{\hrulefill}\thinspace}
\providecommand{\MR}{\relax\ifhmode\unskip\space\fi MR }
\providecommand{\MRhref}[2]{%
  \href{http://www.ams.org/mathscinet-getitem?mr=#1}{#2}
}
\providecommand{\href}[2]{#2}

\end{document}